\documentclass[11pt]{amsart}
\usepackage{amsmath,amssymb,mathrsfs}
\newtheorem{theorem}{Theorem}[section]
\newtheorem{corollary}[theorem]{Corollary}
\newtheorem{lemma}[theorem]{Lemma}
\newtheorem{proposition}[theorem]{Proposition}
\newtheorem*{minoo}{Min-Oo's Conjecture}

\begin{document}

\title{Rigidity phenomena involving scalar curvature}
\author{Simon Brendle}
\address{Department of Mathematics \\ Stanford University \\ 450 Serra Mall, Bldg 380 \\ Stanford, CA 94305} 
\thanks{The author was supported in part by the National Science Foundation under grant DMS-0905628.}
\begin{abstract}
We give a survey of various rigidity results involving scalar curvature. Many of these results are inspired by the positive mass theorem in general relativity. In particular, we discuss the recent solution of Min-Oo's Conjecture for the hemisphere (cf. \cite{Brendle-Marques-Neves}). We also analyze the case of equality in Bray's volume comparison theorem.
\end{abstract}
\maketitle

\section{The positive mass theorem and its geometric consequences}

In this paper, we discuss various rigidity results involving the scalar curvature. Our starting point is the positive mass theorem in general relativity. Recall that a three-manifold $(M,g)$ is said to be asymptotically flat if there exists a compact set $\Omega \subset M$ such that $M \setminus \Omega$ is diffeomorphic to the region $\{x \in \mathbb{R}^3: |x| > 1\}$ and the metric satisfies the estimates 
\begin{align*} 
&|g_{ij}(x) - \delta_{ij}| \leq C \, |x|^{-1}, \\ 
&|\partial_k g_{ij}(x)| \leq C \, |x|^{-2}, \\ 
&|\partial_k \partial_l g_{ij}(x)| \leq C \, |x|^{-3} 
\end{align*} 
for some positive constant $C$. Furthermore, we require that 
\begin{equation}
\label{scalar.curvature.in.L1}
\int_M |R_g| \, d\text{\rm vol}_g < \infty, 
\end{equation}
where $R_g$ denotes the scalar curvature of $(M,g)$. The ADM mass of an asymptotically flat three-manifold $(M,g)$ is defined by 
\begin{equation} 
\label{ADM.mass}
m_{ADM} = \lim_{r \to \infty} \frac{1}{16\pi} \int_{\{|x|=r\}} \sum_{i,j} (\partial_j g_{ij}(x) - \partial_i g_{jj}(x)) \,  \frac{x^i}{r} 
\end{equation} 
(cf. \cite{Arnowitt-Deser-Misner}, \cite{Bartnik1}). It follows from (\ref{scalar.curvature.in.L1}) and the divergence theorem that the limit in (\ref{ADM.mass}) exists, and $m_{ADM}$ is well-defined.

\begin{theorem}[R.~Schoen, S.T.~Yau \cite{Schoen-Yau1}; E.~Witten \cite{Witten}] 
\label{positive.mass.theorem}
Let $(M,g)$ be an asymptotically flat three-manifold with nonnegative scalar curvature. Then the ADM mass of $(M,g)$ is nonnegative. Moreover, if the ADM mass of $(M,g)$ is zero, then $(M,g)$ is isometric to Euclidean space $\mathbb{R}^3$.
\end{theorem}

Theorem \ref{positive.mass.theorem} plays an important role in modern differential geometry. The original proof by Schoen and Yau \cite{Schoen-Yau1} relies on minimal surface techniques. There is an alternative proof due to Witten \cite{Witten}, which uses spinors and the Dirac equation (see also \cite{Parker-Taubes}). 

There is an analogous notion of asymptotic flatness for manifolds of dimension greater than three, and the definition of the ADM mass extends naturally to the higher dimensional setting (cf. \cite{Bartnik1}). The minimal surface arguments of Schoen and Yau imply that the positive mass theorem holds for every asymptotically flat manifold of dimension $n < 8$. On the other hand, by generalizing Witten's argument, Bartnik \cite{Bartnik1} was able to extend the positive mass theorem to spin manifolds of arbitrary dimension. It is an interesting question whether the positive mass theorem holds for non-spin manifolds of dimension $n \geq 8$. This question is studied in recent work of Lohkamp \cite{Lohkamp3}. 

We next consider an important special case of the positive mass theorem. Let $g$ be a Riemannian metric on $\mathbb{R}^n$ which agrees with the Euclidean metric outside a compact set. In this case, $(\mathbb{R}^n,g)$ is asymptotically flat and its ADM mass is equal to zero. Using Bartnik's version of the positive mass theorem, one can draw the following conclusion:

\begin{theorem}
\label{rigidity.of.euclidean.space}
Let $g$ be a metric on $\mathbb{R}^n$ with nonnegative scalar curvature. Moreover, suppose that $g$ agrees with the Euclidean metric outside a compact set. Then $g$ is flat.
\end{theorem}

Similar techniques can be used to show that the $n$-dimensional torus $T^n$ does not admit a metric of positive scalar curvature.

\begin{theorem}[R.~Schoen, S.T.~Yau \cite{Schoen-Yau2}, \cite{Schoen-Yau3}; M.~Gromov, H.B.~Lawson \cite{Gromov-Lawson1}, \cite{Gromov-Lawson2}]
\label{rigidity.of.torus}
Let $g$ be a metric on the torus $T^n$ with nonnegative scalar curvature. Then $g$ is flat.
\end{theorem}

Theorem \ref{rigidity.of.torus} was first proved for $n=3$ by Schoen and Yau \cite{Schoen-Yau2}. The proof relies on minimal surface techniques. In \cite{Schoen-Yau3}, the result was extended to dimension $n < 8$. The general case was settled by Gromov and Lawson using spinor methods (see \cite{Gromov-Lawson1}, \cite{Gromov-Lawson2}).

We next discuss some rigidity results for bounded domains in $\mathbb{R}^n$. It was observed by Miao \cite{Miao} that the positive mass theorem implies the following rigidity result for metrics on the unit ball:

\begin{theorem}
\label{rigidity.for.ball}
Suppose that $g$ is a smooth metric on the unit ball $B^n \subset \mathbb{R}^n$ with the following properties: 
\begin{itemize}
\item The scalar curvature of $g$ is nonnegative.
\item The induced metric on the boundary $\partial B^n$ agrees with the standard metric on $\partial B^n$.
\item The mean curvature of $\partial B^n$ with respect to $g$ is at least $n-1$.
\end{itemize}
Then $g$ is isometric to the standard metric on $B^n$.
\end{theorem}

In 2002, Shi and Tam extended Theorem \ref{rigidity.for.ball} to arbitrary convex domains in $\mathbb{R}^n$. The following result is an important special case of Shi and Tam's theorem (cf. \cite{Shi-Tam}, Theorem 4.1):

\begin{theorem}[Y.~Shi, L.F.~Tam \cite{Shi-Tam}]
\label{shi.tam}
Let $\Omega$ be a strictly convex domain in $\mathbb{R}^n$ with smooth boundary. Moreover, suppose that $g$ is a Riemannian metric on $\Omega$ with the following properties: 
\begin{itemize}
\item The scalar curvature of $g$ is nonnegative.
\item The induced metric on the boundary $\partial \Omega$ agrees with the restriction of the Euclidean metric to $\partial \Omega$.
\item The mean curvature of $\partial \Omega$ with respect to $g$ is positive.
\end{itemize}
Then 
\begin{equation} 
\label{shi.tam.inequality}
\int_{\partial \Omega} (H_0 - H_g) \, d\sigma_g \geq 0, 
\end{equation}
where $H_g$ denotes the mean curvature of $\partial \Omega$ with respect to $g$ and $H_0$ denotes the mean curvature of $\partial \Omega$ with respect to the Euclidean metric. Finally, if equality holds in (\ref{shi.tam.inequality}), then $g$ is flat.
\end{theorem}

In order to prove Theorem \ref{shi.tam}, Shi and Tam glue the metric $g$ to a suitable metric $\tilde{g}$ on the complement $\mathbb{R}^n \setminus \Omega$. The metric $\tilde{g}$ agrees with $g$ along $\partial \Omega$, and the mean curvature of $\partial \Omega$ with respect to $\tilde{g}$ agrees with the mean curvature of $\partial \Omega$ with respect to $g$. Moreover, the metric $\tilde{g}$ is asymptotically flat, and its scalar curvature is equal to zero. Hence, the positive mass theorem implies that $m_{ADM} \geq 0$. On the other hand, Shi and Tam construct a monotone decreasing function $m(r)$ with the property that 
\[m(0) = \int_{\partial \Omega} (H_0 - H_g) \, d\sigma_g\] 
and 
\[\lim_{r \to \infty} m(r) = c(n) \, m_{ADM},\] 
where $c(n)$ is a positive constant (see \cite{Shi-Tam}, Theorem 2.1 and Lemma 4.2). Putting these facts together gives 
\[\int_{\partial \Omega} (H_0 - H_g) \, d\sigma_g \geq c(n) \, m_{ADM} \geq 0,\] 
and the last inequality is strict unless $g$ is flat.

Theorems \ref{rigidity.of.euclidean.space} and \ref{rigidity.of.torus} show that it is not always possible to deform the metric so that the scalar curvature increases at each point. By contrast, Lohkamp \cite{Lohkamp2} proved that for any Riemannian manifold there exist local deformations of the metric which decrease scalar curvature (see also \cite{Lohkamp1}).

\begin{theorem}[J.~Lohkamp \cite{Lohkamp2}]
Let $(M,g)$ be a complete Riemannian manifold of dimension $n \geq 3$, and let $\psi$ be a smooth function on $M$ such that $\psi(x) \leq R_g(x)$ for each point $x \in M$. Let $U = \{x \in M: \psi(x) < R_g(x)\}$, and let $U_\varepsilon$ denote an $\varepsilon$-neighborhood of the set $U$. Given any $\varepsilon > 0$, there exists a smooth metric $\hat{g}$ on $M$ such that $\psi(x) - \varepsilon \leq R_{\hat{g}}(x) \leq \psi(x)$ at each point in $U_\varepsilon$ and $\hat{g} = g$ outside $U_\varepsilon$. 
\end{theorem}

Finally, let us mention a result of Fischer and Marsden \cite{Fischer-Marsden} concerning small deformations of the scalar curvature. To that end, we fix a Riemannian manifold $(M,g)$ of dimension $n \geq 3$. Given any symmetric two-tensor $h$, we define 
\[L_g h = \frac{\partial}{\partial t} R_{g+th} \Big |_{t=0}.\] 
A straightforward calculation gives 
\[L_g h = \sum_{i,j=1}^n (D_{e_i,e_j}^2 h)(e_i,e_j) - \Delta_g(\text{\rm tr}_g(h)) - \langle \text{\rm Ric}_g,h \rangle_g\] 
(see e.g. \cite{Besse},  Theorem 1.174). Consequently, the formal adjoint of $L_g$ is given by 
\[L_g^* f = D^2 f - (\Delta_g f) \, g - f \, \text{\rm Ric}_g.\] 
If the operator $L_g^*$ has trivial kernel, then every function $\psi$ which is sufficiently close to $R_g$ in a suitable sense can be realized as the scalar curvature of a Riemannian metric. A local version of this result was established by Corvino \cite{Corvino}.

We say that $(M,g)$ is static if the operator $L_g^*$ has non-trivial kernel (cf. \cite{Bartnik2}). Examples of static manifolds include the Euclidean space $\mathbb{R}^n$; the hyperbolic space $\mathbb{H}^n$; and the sphere $S^n$ equipped with its standard metric. In the following sections, we will discuss various rigidity theorems for these model spaces.

\section{Rigidity results for hyperbolic space}

In this section, we describe some rigidity results for domains in hyperbolic space. The first result of this type was proved by Min-Oo in 1989. 

\begin{theorem}[M.~Min-Oo \cite{Min-Oo1}]
\label{rigidity.of.hyperbolic.space}
Let $g$ be a metric on $\mathbb{H}^n$ with scalar curvature $R_g \geq -n(n-1)$. Moreover, suppose that $g$ agrees with the hyperbolic metric outside a compact set. Then $g$ has constant sectional curvature $-1$.
\end{theorem}

Theorem \ref{rigidity.of.hyperbolic.space} can be viewed as the analogue of Theorem \ref{rigidity.of.euclidean.space} in the hyperbolic setting. The proof of Theorem \ref{rigidity.of.hyperbolic.space} relies on an adaptation of Witten's proof of the positive mass theorem. The proof does not actually require the metric $g$ to agree with the hyperbolic metric outside a compact set; it suffices to assume that $g$ satisfies certain asymptotic conditions near infinity (see \cite{Min-Oo1} for a precise statement). These asymptotic conditions were later weakened by Andersson and Dahl \cite{Andersson-Dahl}.

We note that there is an analogue of the positive mass theorem for asymptotically hyperbolic manifolds with scalar curvature $R_g \geq -n(n-1)$. Results in this direction were obtained by Chru\'sciel and Herzlich \cite{Chrusciel-Herzlich} and Chru\'sciel and Nagy \cite{Chrusciel-Nagy}. Moreover, the following theorem was established by Wang \cite{Wang}:

\begin{theorem}[X.~Wang \cite{Wang}]
Let $g$ be a metric on the unit ball $B^n$ with scalar curvature $R_g \geq -n(n-1)$. Moreover, suppose that $g$ satisfies an asymptotic expansion of the form 
\[g = \sinh^{-2}(r) \, \Big ( dr^2 + g_0 + \frac{r^n}{n} \, h + O(r^{n+1}) \Big ),\] 
where $g_0$ denotes the round metric on the boundary $\partial B^n = S^{n-1}$ and $r$ is a boundary defining function. Then 
\[\int_{S^{n-1}} \text{\rm tr}_{g_0}(h) \, d\text{\rm vol}_{g_0} \geq \bigg | \int_{S^{n-1}} \text{\rm tr}_{g_0}(h) \, x \, d\text{\rm vol}_{g_0} \bigg |.\] 
Moreover, if equality holds, then $g$ is isometric to the hyperbolic metric.
\end{theorem}

Finally, we point out that Boualem and Herzlich have obtained similar rigidity results for K\"ahler manifolds that are asymptotic to complex hyperbolic space in a suitable sense (cf. \cite{Boualem-Herzlich}, \cite{Herzlich}).

\section{Min-Oo's Conjecture for the hemisphere}

In this section, we discuss rigidity questions for the hemisphere $S^n$. For abbreviation, we will denote by $\overline{g}$ the standard metric on $S^n$. Motivated by the positive mass theorem and its analogue in the asymptotically hyperbolic setting, Min-Oo proposed the following conjecture (cf. \cite{Min-Oo2}, Theorem 4):

\begin{minoo}
Suppose that $g$ is a smooth metric on the hemisphere $S_+^n = \{x \in S^n: x_{n+1} \geq 0\}$ with the following properties: 
\begin{itemize}
\item The scalar curvature of $g$ is at least $n(n-1)$.
\item The induced metric on the boundary $\partial S_+^n$ agrees with the standard metric on $\partial S_+^n$.
\item The boundary $\partial S_+^n$ is totally geodesic with respect to $g$.
\end{itemize}
Then $g$ is isometric to the standard metric on $S_+^n$.
\end{minoo}

Min-Oo's Conjecture is very natural given the analogy with the positive mass theorem, and was widely expected to be true; see e.g. \cite{Gromov}, p.~47, or \cite{Hang-Wang2}, p.~629. Various attempts have been made to prove it (using both spinor and minimal surface techniques), and many partial results have been obtained. In particular, it follows from a classical result of Toponogov that Min-Oo's Conjecture holds in dimension $2$.

\begin{theorem}[V.~Toponogov \cite{Toponogov}]
\label{dim.2}
Let $(M,g)$ be a compact surface with totally geodesic boundary $\partial M$. If the Gaussian curvature of $(M,g)$ satisfies $K \geq 1$, then the length of $\partial M$ is at most $2\pi$. Moreover, if equality holds, then $(M,g)$ is isometric to the hemisphere $S_+^2$ equipped with its standard metric.
\end{theorem}

An alternative proof of Theorem \ref{dim.2} was given by Hang and Wang (cf. \cite{Hang-Wang2}, Theorem 4). In higher dimensions, Hang and Wang \cite{Hang-Wang2} showed that Min-Oo's Conjecture holds if the lower bound for the scalar curvature is replaced by a lower bound for the Ricci tensor.

\begin{theorem}[F.~Hang, X.~Wang \cite{Hang-Wang2}]
\label{hangwang.thm}
Suppose that $g$ is a smooth metric on the hemisphere $S_+^n$ with the following properties: 
\begin{itemize}
\item The Ricci curvature of $g$ is bounded from below by $\text{\rm Ric}_g \geq (n-1) \, g$.
\item The induced metric on the boundary $\partial S_+^n$ agrees with the standard metric on $\partial S_+^n$.
\item The second fundamental form of the boundary $\partial S_+^n$ with respect to $g$ is nonnegative.
\end{itemize}
Then $g$ is isometric to the standard metric on $S_+^n$.
\end{theorem}

The proof of Theorem \ref{hangwang.thm} relies on an interesting application of Reilly's formula. Moreover, Hang and Wang \cite{Hang-Wang1} were able to verify Min-Oo's Conjecture for metrics conformal to the standard metric. In particular, Min-Oo's Conjecture is true if the metric $g$ is rotationally symmetric.

\begin{theorem}[F.~Hang, X.~Wang \cite{Hang-Wang1}]
\label{conf.flat}
Suppose that $g = e^{2w} \, \overline{g}$ is a metric in the conformal class of $\overline{g}$ with scalar curvature $R_g \geq n(n-1)$. If $g = \overline{g}$ along the boundary $\partial S_+^n$, then $g = \overline{g}$ at each point in $S_+^n$.
\end{theorem}

Note that Theorem \ref{conf.flat} does not require any assumptions on the second fundamental form of $\partial S_+^n$.

In a joint work with F.C. Marques, we have obtained scalar curvature rigidity results for certain geodesic balls in $S^n$; see \cite{Brendle-Marques2}.

\begin{theorem}[S.~Brendle, F.C.~Marques \cite{Brendle-Marques2}]
\label{geodesic.balls}
Fix a real number $c \geq \frac{2}{\sqrt{n+3}}$, and let $\Omega = \{x \in S^n: x_{n+1} \geq c\}$. Moreover, suppose that $g$ is a Riemannian metric on $\Omega$ with the following properties: 
\begin{itemize}
\item $R_g \geq n(n-1)$ at each point in $\Omega$. 
\item The metrics $g$ and $\overline{g}$ induce the same metric on $\partial \Omega$. 
\item $H_g \geq H_{\overline{g}}$ at each point on $\partial \Omega$. 
\end{itemize} 
If $g - \overline{g}$ is sufficiently small in the $C^2$-norm, then $\varphi^*(g) = \overline{g}$ for some diffeomorphism $\varphi: \Omega \to \Omega$ with $\varphi|_{\partial \Omega} = \text{\rm id}$.
\end{theorem}

The following result is an immediate consequence of Theorem \ref{geodesic.balls}.

\begin{corollary}[S.~Brendle, F.C.~Marques \cite{Brendle-Marques2}]
Suppose that $g$ is a Riemannian metric on the hemisphere $S_+^n$ with the following properties: 
\begin{itemize}
\item $R_g \geq n(n-1)$ at each point in $S_+^n$. 
\item The metrics $g$ and $\overline{g}$ agree in the region $\big \{ x \in S^n: 0 \leq x_{n+1} \leq \frac{2}{\sqrt{n+3}} \big \}$. 
\end{itemize} 
If $g - \overline{g}$ is sufficiently small in the $C^2$-norm, then $g$ is isometric to the standard metric $\overline{g}$.
\end{corollary}

Eichmair \cite{Eichmair} has verified Min-Oo's Conjecture for $n=3$, assuming that the boundary satisfies an isoperimetric condition. The proof of this theorem uses techniques developed by Bray \cite{Bray}. Moreover, Huang and Wu \cite{Huang-Wu} showed that Min-Oo's Conjecture holds for graphs in $\mathbb{R}^{n+1}$.

We next describe a non-rigidity theorem for the hemisphere. This result implies that Min-Oo's Conjecture fails in dimension $n \geq 3$. The argument involves two steps. In a first step, we perturb the standard metric on the hemisphere $S_+^n$ in such a way that the scalar curvature increases at each point and the mean curvature of the boundary becomes positive:

\begin{theorem}[S.~Brendle, F.C.~Marques, A.~Neves \cite{Brendle-Marques-Neves}]
\label{counterexample.A}
Given any integer $n \geq 3$, there exists a smooth metric $g$ on the hemisphere $S_+^n$ with the following properties:
\begin{itemize}
\item The scalar curvature of $g$ is strictly greater than $n(n-1)$.
\item At each point on $\partial S_+^n$, we have $g - \overline{g} = 0$, where $\overline{g}$ denotes the standard metric on $S_+^n$.
\item The mean curvature of $\partial S_+^n$ with respect to $g$ is strictly positive (i.e. the mean curvature vector is inward-pointing).
\end{itemize}
\end{theorem}

The proof of Theorem \ref{counterexample.A} relies on a perturbation analysis. This construction is inspired by the counterexamples to Schoen's Compactness Conjecture for the Yamabe problem (cf. \cite{Brendle}, \cite{Brendle-Marques1}). We will give an outline of the proof of Theorem \ref{counterexample.A} in Section \ref{construction.of.counterexamples}.

In a second step, we construct metrics on the hemisphere $S_+^n$ which have scalar curvature at least $n(n-1)$ and agree with the standard metric in a neighborhood of the boundary. To that end, we glue the metrics constructed in Theorem \ref{counterexample.A} to a rotationally symmetric model metric. The proof relies on a general gluing theorem, which is of interest in itself.

\begin{theorem}[S.~Brendle, F.C.~Marques, A.~Neves \cite{Brendle-Marques-Neves}]
\label{counterexample.B}
Let $M$ be a compact manifold of dimension $n$ with boundary $\partial M$, and let $g$ and $\tilde{g}$ be two smooth Riemannian metrics on $M$ such that $g - \tilde{g} = 0$ at each point on $\partial M$. Moreover, we assume that $H_g - H_{\tilde{g}} > 0$ at each point on $\partial M$. Given any real number $\varepsilon > 0$ and any neighborhood $U$ of $\partial M$, there exists a smooth metric $\hat{g}$ on $M$ with the following properties: 
\begin{itemize}
\item We have the pointwise inequality $R_{\hat{g}}(x) \geq \min \{R_g(x),R_{\tilde{g}}(x)\} - \varepsilon$ at each point $x \in M$.
\item $\hat{g}$ agrees with $g$ outside $U$.
\item $\hat{g}$ agrees with $\tilde{g}$ in a neighborhood of $\partial M$.
\end{itemize}
\end{theorem}

The proof of Theorem \ref{counterexample.B} involves a delicate choice of cut-off functions (see \cite{Brendle-Marques-Neves}, Section 4). We will omit the details here.

Combining Theorem \ref{counterexample.A} and Theorem \ref{counterexample.B}, we are able to construct counterexamples to Min-Oo's Conjecture in dimension $n \geq 3$.

\begin{theorem}[S.~Brendle, F.C.~Marques, A.~Neves \cite{Brendle-Marques-Neves}]
\label{counterexample.C}
Given any integer $n \geq 3$, there exists a smooth metric $\hat{g}$ on the hemisphere $S_+^n$ with the following properties:
\begin{itemize}
\item The scalar curvature of $\hat{g}$ is at least $n(n-1)$ at each point on $S_+^n$.
\item The scalar curvature of $\hat{g}$ is strictly greater than $n(n-1)$ at some point on $S_+^n$.
\item The metric $\hat{g}$ agrees with the standard metric $\overline{g}$ in a neighborhood of $\partial S_+^n$.
\end{itemize} 
\end{theorem}

\textit{Sketch of the proof of Theorem \ref{counterexample.C}.} 
Let $\delta$ be a small positive number. We first construct a rotationally symmetric metric $\tilde{g}_\delta$ on the hemisphere $S_+^n$ such that $R_{\tilde{g}_\delta} > n(n-1)$ in the region $\{x \in S^n: \delta < x_{n+1} < 3\delta\}$ and $\tilde{g}_\delta = \overline{g}$ in the region $\{x \in S^n: 0 \leq x_{n+1} \leq \delta\}$. 

Let $M_\delta = \{x \in S^n: x_{n+1} \geq 2\delta\}$. Using Theorem \ref{counterexample.B}, we can construct a metric $g_\delta$ on $M_\delta$ with the following properties: 
\begin{itemize}
\item $R_{g_\delta} > n(n-1)$ at each point in $M_\delta$.
\item $g_\delta - \tilde{g}_\delta = 0$ at each point on the boundary $\partial M_\delta$. 
\item $H_{g_\delta} - H_{\tilde{g}_\delta} > 0$ at each point on $\partial M_\delta$.
\end{itemize}
Applying Theorem \ref{counterexample.B} to the metrics $g_\delta$ and $\tilde{g}_\delta$, we obtain a metric $\hat{g}$ on $M_\delta$ with the property that $R_{\hat{g}} > n(n-1)$ at each point in $M_\delta$ and $\hat{g} = \tilde{g}_\delta$ in a  neighborhood of $\partial M_\delta$. Hence, we may extend the metric $\hat{g}$ to the hemisphere in such a way that the resulting metric has all the required properties. From this, Theorem \ref{counterexample.C} follows. \\

To conclude this section, we state a corollary of Theorem \ref{counterexample.C}. Let $\hat{g}$ be the metric constructed in Theorem \ref{counterexample.C}. We may extend $\hat{g}$ to a metric on $S^n$ which is invariant under antipodal reflection. The resulting metric descends to a metric on the real projective space $\mathbb{RP}^n$. Hence, we can draw the following conclusion:

\begin{corollary}[S.~Brendle, F.C.~Marques, A.~Neves \cite{Brendle-Marques-Neves}]
\label{projective.space}
Given any integer $n \geq 3$, there exists a smooth metric $g$ on the real projective space $\mathbb{RP}^n$ with the following properties:
\begin{itemize}
\item The scalar curvature of $g$ is at least $n(n-1)$ at each point on $\mathbb{RP}^n$.
\item The scalar curvature of $g$ is strictly greater than $n(n-1)$ at some point on $\mathbb{RP}^n$.
\item The metric $g$ agrees with the standard metric in a neighborhood of the equator in $\mathbb{RP}^n$.
\end{itemize}
\end{corollary}

\section{Sketch of the proof of Theorem \ref{counterexample.A}}

\label{construction.of.counterexamples}

In this section, we describe the main ideas involved in the proof of Theorem \ref{counterexample.A}. A complete proof is presented in \cite{Brendle-Marques-Neves}. A crucial issue is that the standard metric $\overline{g}$ on $S^n$ is static. In fact, if we denote by $f: S^n \to \mathbb{R}$ the restriction of the coordinate function $x_{n+1}$ to the unit sphere $S^n$, then $f$ satisfies the equation 
\begin{equation} 
\label{static}
\overline{D}^2 f - (\Delta_{\overline{g}} f) \, \overline{g} - f \, \text{\rm Ric}_{\overline{g}} = 0. 
\end{equation} 
Consequently, the linearized operator $L_{\overline{g}}$ fails to be surjective. In particular, Corvino's theorem concerning local deformations of the scalar curvature (see \cite{Corvino}, Theorem 1) does not apply in this situation.

We next consider the upper hemisphere $S_+^n = \{f \geq 0\}$. To fix notation, let $\Sigma = \{f = 0\}$ denote the equator in $S^n$ and let $\nu$ be the outward-pointing unit normal vector field along $\Sigma$. Given any Riemannian metric $g$ on $S_+^n$, we define 
\begin{equation} 
\mathscr{F}(g) = \int_{S_+^n} R_g \, f \, d\text{\rm vol}_{\overline{g}} + 2 \, \text{\rm area}(\Sigma,g). 
\end{equation}
This defines a functional $\mathscr{F}$ on the space of Riemannian metrics on $S_+^n$. Using the relation (\ref{static}), one can show that the first variation of $\mathscr{F}$ at $\overline{g}$ vanishes. More precisely, we have the following result: 

\begin{proposition}[\cite{Brendle-Marques-Neves}, Proposition 9]
\label{first.variation}
Let $g(t)$ be a smooth one-parameter family of Riemannian metrics on $S_+^n$ with $g(0) = \overline{g}$. Then $\frac{d}{dt} \mathscr{F}(g(t)) \big |_{t=0} = 0$.
\end{proposition}

Note that Proposition \ref{first.variation} holds for arbitrary variations of the metric, including those that change the induced metric on the boundary. This fact will play a crucial role in the argument.

The strategy is to deform the standard metric $\overline{g}$ in such a way that the scalar curvature is unchanged to first order and the second variation of the functional $\mathscr{F}$ is positive. In order to construct variations with this property, we need the following auxiliary result:

\begin{proposition}[\cite{Brendle-Marques-Neves}, Proposition 10]
\label{eta}
Assume that $n \geq 3$. Then there exists a function $\eta: \Sigma \to \mathbb{R}$ such that 
\[\Delta_\Sigma \eta + (n-1) \eta < 0\] 
and 
\[\int_\Sigma (|\nabla_\Sigma \eta|^2 - (n-1) \eta^2) \, d\sigma_{\overline{g}} > 0.\] 
\end{proposition}

\textit{Sketch of the proof of Proposition \ref{eta}.} 
We define a function $\eta: \Sigma \to \mathbb{R}$ by 
\[\eta = -c - 1 + \frac{n-1}{2} \, x_n^2 + \frac{(n-1)(n+1)}{24} \, x_n^4 + \frac{(n-1)(n+1)(n+3)}{240} \, x_n^6,\] 
where $c$ is a small positive constant. Then 
\[\Delta_\Sigma \eta + (n-1) \eta = -(n-1)c - \frac{(n-1)(n+1)(n+3)(n+5)}{48} \, x_n^6 < 0.\] 
Moreover, a straightforward calculation shows that  
\[\int_\Sigma (|\nabla_\Sigma \eta|^2 - (n-1) \eta^2) \, d\sigma_{\overline{g}} > 0\] 
if $c > 0$ is sufficiently small. Hence, the function $\eta$ has the required properties. This completes the proof of Proposition \ref{eta}. \\

Note that the proof of Proposition \ref{eta} fails for $n = 2$. Proposition \ref{eta} has a natural geometric interpretation: it implies that, for $n \geq 3$, there exist deformations of the equator in $S^n$ which increase area and have positive mean curvature.

\begin{corollary}
Assume that $n \geq 3$. Then there exists a one-parameter family of hypersurfaces $\Sigma_t$ with the following properties: 
\begin{itemize} 
\item $\Sigma_0 = \Sigma$.
\item $\Sigma_t$ has positive mean curvature for each $t > 0$.
\item $\frac{d^2}{dt^2} \text{\rm area}(\Sigma_t,\overline{g}) \big |_{t=0} > 0$.
\end{itemize}
\end{corollary}

In the remainder of this section, we will always assume that $n \geq 3$. Let $\eta: \Sigma \to \mathbb{R}$ be the function constructed in Proposition \ref{eta}. We can find a smooth vector field $X$ on $S^n$ such that 
\[X = \eta \, \nu\] 
and 
\[\mathscr{L}_X \overline{g} = 0\] 
at each point on $\Sigma$. The exact choice of $X$ is not important; all that matters is the behavior of $X$ near the equator.

The vector field $X$ generates a one-parameter group of diffeomorphisms, which we denote by $\varphi_t: S^n \to S^n$. For each $t$, we define two Riemannian metrics $g_0(t)$ and $g_1(t)$ by 
\[g_0(t) = \overline{g} + t \, \mathscr{L}_X \overline{g}\] 
and 
\[g_1(t) = \varphi_t^*(\overline{g}).\] 
It follows from our choice of $X$ that the metric $g_0(t)$ agrees with the standard metric $\overline{g}$ at each point on $\Sigma$. By contrast, the metric $g_1(t)$ does not agree with $\overline{g}$ along $\Sigma$.

\begin{proposition}[\cite{Brendle-Marques-Neves}, Proposition 11]
\label{key}
There exists a smooth function $Q$ on the hemisphere $S_+^n$ such that 
\[R_{g_0(t)} = n(n-1) + \frac{1}{2} \, t^2 \, Q + O(t^3)\] 
and 
\[\int_{S_+^n} Q \, f \, d\text{\rm vol}_{\overline{g}} > 0.\] 
\end{proposition}

\textit{Sketch of the proof of Proposition \ref{key}.} 
Note that $g_0(t) = g_1(t) + O(t^2)$. Since $R_{g_1(t)} = n(n-1)$, it follows that $R_{g_0(t)} = n(n-1) + O(t^2)$. For abbreviation, let 
\[Q = \frac{\partial^2}{\partial t^2} R_{g_0(t)} \Big |_{t=0}.\] 
By Proposition \ref{first.variation}, the first variation of the functional $\mathscr{F}$ at $\overline{g}$ vanishes. This implies 
\begin{equation} 
\label{key.identity}
\frac{d^2}{dt^2} \mathscr{F}(g_0(t)) \Big |_{t=0} = \frac{d^2}{dt^2} \mathscr{F}(g_1(t)) \Big |_{t=0}. 
\end{equation} 
The left hand side in (\ref{key.identity}) is given by 
\[\frac{d^2}{dt^2} \mathscr{F}(g_0(t)) \Big |_{t=0} =  \int_{S_+^n} Q \, f \, d\text{\rm vol}_{\overline{g}}.\] 
On the other hand, using the identity $R_{g_1(t)} = n(n-1)$, the right hand side in (\ref{key.identity}) can be rewritten as 
\begin{align*} 
\frac{d^2}{dt^2} \mathscr{F}(g_1(t)) \Big |_{t=0} 
&= 2 \, \frac{d^2}{dt^2} \text{\rm area}(\Sigma,g_1(t)) \Big |_{t=0} \\ 
&= 2 \, \frac{d^2}{dt^2} \text{\rm area}(\varphi_t(\Sigma),\overline{g}) \Big |_{t=0} \\ 
&= 2 \int_\Sigma (|\nabla_\Sigma \eta|^2 - (n-1) \eta^2) \, d\sigma_{\overline{g}}. 
\end{align*} 
Putting these facts together, we obtain 
\[ \int_{S_+^n} Q \, f \, d\text{\rm vol}_{\overline{g}} = 2 \int_\Sigma (|\nabla_\Sigma \eta|^2 - (n-1) \eta^2) \, d\sigma_{\overline{g}} > 0,\] 
completing the proof of Proposition \ref{key}. \\

We next consider the elliptic equation 
\begin{equation} 
\Delta_{\overline{g}} u + nu = Q - \frac{\int_{S_+^n} Q \, f \, d\text{\rm vol}_{\overline{g}}}{\int_{S_+^n} f \, d\text{\rm vol}_{\overline{g}}} 
\end{equation}
with Dirichlet boundary condition $u|_\Sigma = 0$. This boundary value problem has a smooth solution $u: S_+^n \to \mathbb{R}$. We now define 
\[g(t) = \overline{g} + t \, \mathscr{L}_X \overline{g} + \frac{1}{2(n-1)} \, t^2 \, u \, \overline{g}.\] 
The metric $g(t)$ agrees with the standard metric $\overline{g}$ at each point on $\Sigma$. Moreover, the scalar curvature of $g(t)$ is given by 
\begin{align*} 
R_{g(t)} 
&= R_{g_0(t)} - \frac{1}{2} \, t^2 \, (\Delta_{\overline{g}} u + nu) + O(t^3) \\ 
&= n(n-1) + \frac{1}{2} \, t^2 \, Q  - \frac{1}{2} \, t^2 \, (\Delta_{\overline{g}} u + nu) + O(t^3) \\ 
&= n(n-1) + \frac{1}{2} \, t^2 \, \frac{\int_{S_+^n} Q \, f \, d\text{\rm vol}_{\overline{g}}}{\int_{S_+^n} f \, d\text{\rm vol}_{\overline{g}}} + O(t^3). 
\end{align*} 
Since $\int_{S_+^n} Q \, f \, d\text{\rm vol}_{\overline{g}} > 0$, it follows that 
\[\inf_{S_+^n} R_{g(t)} > n(n-1)\] 
if $t > 0$ is sufficiently small.

It remains to compute the mean curvature of $\Sigma$ with respect to the metric $g(t)$. Since $g(t) = g_1(t) + O(t^2)$, we have 
\[H_{g(t)} = H_{g_1(t)} + O(t^2) =  -t \, (\Delta_\Sigma \eta + (n-1)\eta) + O(t^2).\] 
By Proposition \ref{eta}, we have $\Delta_\Sigma \eta + (n-1)\eta < 0$ at each point on $\Sigma$. Consequently, we have 
\[\inf_\Sigma H_{g(t)} > 0\] 
if $t > 0$ is sufficiently small. 

\section{Other rigidity results involving scalar curvature}

In this section, we describe other related rigidity results. 

The following result was established by Llarull \cite{Llarull} (see also the survey paper \cite{Gromov}). 

\begin{theorem}[M.~Llarull \cite{Llarull}]
\label{Llarull.thm}
Let $g$ be a Riemannian metric on $S^n$ with scalar curvature $R_g \geq n(n-1)$. Moreover, suppose that $g \geq \overline{g}$ at each point on $S^n$. Then $g = \overline{g}$ at each point on $S^n$.
\end{theorem}

In the even-dimensional case, Listing \cite{Listing} was able to generalize Theorem \ref{Llarull.thm} as follows: 

\begin{theorem}[M.~Listing \cite{Listing}]
Let $n \geq 4$ be an even integer. Moreover, suppose that $g$ is a Riemannian metric on $S^n$ satisfying $R_g \geq (n-1) \, \text{\rm tr}_g(\overline{g})$ at each point on $S^n$. Then $g$ is a constant multiple of the standard metric $\overline{g}$.
\end{theorem}

The following theorem due to Bray gives a sharp upper bound for the volume of a three-manifold with scalar curvature at least $6$.

\begin{theorem}[H.~Bray \cite{Bray}] 
\label{vol.comp}
Fix a real number $\varepsilon \in (0,1)$ with the property that 
\begin{align} 
\label{varepsilon}
&\int_0^y \big ( 36\pi - 27(1-\varepsilon) \, y^{\frac{2}{3}} - 9\varepsilon \, x^{\frac{2}{3}} \big )^{-\frac{1}{2}} \, dx \notag \\ 
&+ \int_y^{z^{\frac{3}{2}}} \big ( 36\pi - 18(1-\varepsilon) \, y \, x^{-\frac{1}{3}} - 9 \, x^{\frac{2}{3}} \big )^{-\frac{1}{2}} \, dx < \pi^2 
\end{align} 
for all pairs $(y,z) \in \mathbb{R} \times \big [ \frac{4\pi}{3-2\varepsilon},4\pi \big )$ satisfying $2(1-\varepsilon) \, y = z^{\frac{1}{2}} \, (4\pi - z)$. Moreover, let $(M,g)$ be a compact three-manifold satisfying $R_g \geq 6$ and $\text{\rm Ric}_g \geq 2\varepsilon \, g$. If $\text{\rm vol}(M,g) \geq \text{\rm vol}(S^3,\overline{g})$, then $(M,g)$ is isometric to $(S^3,\overline{g})$.
\end{theorem}

Theorem \ref{vol.comp} is slightly stronger than the result proved in \cite{Bray}. It follows from Theorem 19 in \cite{Bray} that every three-manifold $(M,g)$ with $R_g \geq 6$ and $\text{\rm Ric}_g \geq 2\varepsilon \, g$ satisfies $\text{\rm vol}(M,g) \leq \text{\rm vol}(S^3,\overline{g})$. However, the case of equality has not been discussed in the literature. For the convenience of the reader, we shall provide a detailed exposition in the following section.

Gursky and Viaclovsky proved that the condition (\ref{varepsilon}) is satisfied for $\varepsilon = \frac{1}{2}$ (see \cite{Gursky-Viaclovsky}, Section 4.1). Hence, we can draw the following conclusion:

\begin{corollary}
Let $(M,g)$ be a compact three-manifold satisfying $R_g \geq 6$ and $\text{\rm Ric}_g \geq g$. If $\text{\rm vol}(M,g) \geq \text{\rm vol}(S^3,\overline{g})$, then $(M,g)$ is isometric to $(S^3,\overline{g})$.
\end{corollary}

We now discuss a rigidity result for Riemannian metrics on $\mathbb{RP}^3$ with scalar curvature at least $6$. To fix notation, we denote by $\mathscr{F}$ the set of all embedded surfaces $\Sigma \subset \mathbb{RP}^3$ with the property that $\Sigma$ is homeomorphic to $\mathbb{RP}^2$.

\begin{theorem}[H.~Bray, S.~Brendle, M.~Eichmair, A.~Neves \cite{Bray-Brendle-Eichmair-Neves}]
\label{rp2}
Let $g$ be a Riemannian metric on $\mathbb{RP}^3$ with scalar curvature $R_g \geq 6$. Moreover, suppose that $\Sigma \in \mathscr{F}$ is a surface which has minimal area among all surfaces in $\mathscr{F}$. Then $\text{\rm area}(\Sigma,g) \leq 2\pi$. Moreover, if equality holds, then $g$ is isometric to the standard metric on $\mathbb{RP}^3$.
\end{theorem}

\textit{Sketch of the proof of Theorem \ref{rp2}.} Given any metric $g$ on $\mathbb{RP}^3$ and any stable minimal surface $\Sigma \in \mathscr{F}$, one can show that 
\begin{equation} 
\label{main.inequality}
\text{\rm area}(\Sigma,g) \, \inf_{\mathbb{RP}^3} R_g \leq 12\pi 
\end{equation} 
(cf. \cite{Bray-Brendle-Eichmair-Neves}, Corollary 8). To prove (\ref{main.inequality}), we use special choices of variations in the second variation formula. These variations are constructed by adapting a technique of Hersch \cite{Hersch}.

We now sketch the proof of the rigidity statement. Suppose that $g$ is a metric on $\mathbb{RP}^3$ with scalar curvature $R_g \geq 6$, and $\Sigma \in \mathscr{F}$ is an area-minizing surface with $\text{\rm area}(\Sigma,g) = 2\pi$. Let $\tilde{g}(t)$, $t \in [0,T)$, denote the unique solution to the Ricci flow with initial metric $\tilde{g}(0) = g$. Fix a real number $\tau \in (0,T)$. By a theorem of Meeks, Simon, and Yau \cite{Meeks-Simon-Yau}, we can find a surface $\tilde{\Sigma} \in \mathscr{F}$ which has minimal area with respect to the metric $\tilde{g}(\tau)$. The key idea is to establish a lower bound for the area of $\tilde{\Sigma}$. More precisely, we have 
\begin{equation} 
\label{ineq.1}
\text{\rm area}(\tilde{\Sigma},\tilde{g}(\tau)) \geq \text{\rm area}(\Sigma,g) - 8\pi \tau = 2\pi \, (1 - 4\tau) 
\end{equation}
(cf. \cite{Bray-Brendle-Eichmair-Neves}, Proposition 10). On the other hand, it follows from the maximum principle that 
\begin{equation} 
\label{ineq.2}
\inf_{\mathbb{RP}^3} R_{\tilde{g}(\tau)} \geq \frac{6}{1-4\tau}. 
\end{equation} 
Moreover, applying (\ref{main.inequality}) to the metric $\tilde{g}(\tau)$ gives 
\[\text{\rm area}(\tilde{\Sigma},\tilde{g}(\tau)) \, \inf_{\mathbb{RP}^3} R_{\tilde{g}(\tau)} \leq 12\pi.\] 
Therefore, the inequalities (\ref{ineq.1}) and (\ref{ineq.2}) hold as equalities. The strict maximum principle then implies that $g$ has constant sectional curvature $1$. This completes the proof of Theorem \ref{rp2}. \\

By Corollary \ref{projective.space}, there are non-trivial examples of Riemannian metrics on $\mathbb{RP}^3$ that have scalar curvature at least $6$ and agree with the standard metric in a neighborhood of the equator. In this case, the equator is a stable minimal surface of area $2\pi$. Therefore, the rigidity statement in Theorem \ref{rp2} is no longer true if we replace the assumption that $\Sigma$ is area-minimizing by the weaker condition that $\Sigma$ is stable. 

We next describe an analogous estimate for the area of area-minimizing two-spheres in three-manifolds. 

\begin{theorem}[H.~Bray, S.~Brendle, A.~Neves \cite{Bray-Brendle-Neves}]
\label{cylinder}
Let $(M,g)$ be a compact three-manifold with scalar curvature $R_g \geq 2$. Moreover, suppose that $\Sigma$ is an immersed two-sphere which minimizes area in its homotopy class. Then $\text{\rm area}(\Sigma,g) \leq 4\pi$. Moreover, if equality holds, then the universal cover of $(M,g)$ is isometric to the cylinder $S^2 \times \mathbb{R}$ equipped with its standard metric.
\end{theorem}

\textit{Sketch of the proof of Theorem \ref{cylinder}.} 
The inequality $\text{\rm area}(\Sigma,g) \leq 4\pi$ follows from the stability inequality. We now describe of the rigidity statement. Let $\Sigma$ be an area-minimizing two-sphere with $\text{\rm area}(\Sigma,g) = 4\pi$, and let $\nu$ be a unit normal vector field along $\Sigma$. It is easy to see that $\Sigma$ is totally geodesic and $\text{\rm Ric}(\nu,\nu) = 0$ at each point on $\Sigma$. Using the implicit function theorem, we may construct a one-parameter family of stable constant mean curvature surfaces 
\[\Sigma_t = \{\exp_x(w(x,t) \, \nu(x)): x \in \Sigma\},\] 
where $w(x,0) = 0$ and $\frac{\partial}{\partial t} w(x,t) \big |_{t=0} = 1$. Let $\nu_t$ denote the unit normal vector to $\Sigma_t$. We assume that $\nu_t$ is chosen so that $\nu_0 = \nu$. The mean curvature vector of $\Sigma_t$ can be written as $-H(t) \, \nu_t$, where $H(t)$ depends only on $t$. By definition of $\Sigma$, we have 
\[\text{\rm area}(\Sigma_t,g) \geq \text{\rm area}(\Sigma,g) = 4\pi.\] 
Using the inequality $R \geq 2$ and the Gauss equations, we obtain 
\[\int_{\Sigma_t} (\text{\rm Ric}(\nu_t,\nu_t) + |I\!I_t|^2) \geq 0.\] 
Combining this relation with the stability inequality, one can show that $H'(t) \leq 0$ for all $t \in (-\delta,\delta)$ (see \cite{Bray-Brendle-Neves}, p.~827). Since $H(0) = 0$, we conclude that 
\[\text{\rm area}(\Sigma_t,g) \leq \text{\rm area}(\Sigma,g) = 4\pi\] 
for $t \in (-\delta,\delta)$. Consequently, $\Sigma_t$ must be totally geodesic, and we have $\text{\rm Ric}(\nu_t,\nu_t) = 0$ at each point on $\Sigma_t$. From this, we deduce that $(M,g)$ locally splits as a product. \\

As above, the rigidity statement in Theorem \ref{cylinder} fails if we replace the condition that $\Sigma$ is area-minimizing by the weaker condition that $\Sigma$ is stable. Moreover, the proof of Theorem \ref{cylinder} can be adapted to other settings, see e.g. \cite{Micallef-Moraru}, \cite{Nunes}. 

Finally, we note that Cai and Galloway \cite{Cai-Galloway} have obtained an interesting rigidity result for minimal tori in three-manifolds with nonnegative scalar curvature (see also \cite{Fischer-Colbrie-Schoen}): 

\begin{theorem}[M.~Cai, G.~Galloway \cite{Cai-Galloway}]
Let $(M,g)$ be a three-manifold with nonnegative scalar curvature, and let $\Sigma$ be a two-sided minimal torus in $(M,g)$ which is locally area-minimizing. Then $g$ is flat in a neighborhood of $\Sigma$.
\end{theorem}

\section{The case of equality in Bray's volume comparison theorem}

In this final section, we describe the proof of Theorem \ref{vol.comp}. We will closely follow the original argument of Bray \cite{Bray}. Let us fix a real number $\varepsilon \in (0,1)$ satisfying (\ref{varepsilon}). For each $z \in \big ( \frac{4\pi}{3-2\varepsilon},4\pi \big ]$, we define 
\[y(z) = \frac{1}{2(1-\varepsilon)} \, z^{\frac{1}{2}} \, (4\pi - z) \in [0,z^{\frac{3}{2}})\] 
and 
\begin{align*} 
\gamma(z) 
&= \int_0^{y(z)} \big ( 36\pi - 27(1-\varepsilon) \, y(z)^{\frac{2}{3}} - 9\varepsilon \, x^{\frac{2}{3}} \big )^{-\frac{1}{2}} \, dx \\ 
&+ \int_{y(z)}^{z^{\frac{3}{2}}} \big ( 36\pi - 18(1-\varepsilon) \, y(z) \, x^{-\frac{1}{3}} - 9 \, x^{\frac{2}{3}} \big )^{-\frac{1}{2}} \, dx. 
\end{align*} 
Moreover, for each $z \in \big ( 0,\frac{4\pi}{3-2\varepsilon} \big ]$, we define 
\[\gamma(z) = \int_0^{z^{\frac{3}{2}}} \big ( 9\varepsilon \, z - 9\varepsilon \, x^{\frac{2}{3}} \big )^{-\frac{1}{2}} \, dx.\] 
Note that the function $\gamma(z)$ is continuous.

\begin{proposition}
\label{comparison.function} 
For each $z \in (0,4\pi]$, there exists a function $f_z: (-\gamma(z),\gamma(z)) \to (0,z^{\frac{3}{2}}]$ with the following properties: 
\begin{itemize}
\item The function $f_z$ is twice continuously differentiable.
\item For each $s \in (-\gamma(z),\gamma(z))$, we have 
\begin{equation} 
\label{ode.for.f}
f_z''(s) = \min \Big \{ \frac{36\pi - f_z'(s)^2}{6f_z(s)} - \frac{9}{2} \, f_z(s)^{-\frac{1}{3}},-3\varepsilon \, f_z(s)^{-\frac{1}{3}} \Big \} 
\end{equation}
\item $f_z(0) = z^{\frac{3}{2}}$ and $f_z'(0) = 0$. 
\item $\lim_{s \to \gamma(z)} f_z(s) = \lim_{s \to -\gamma(z)} f_z(s) = 0$.
\end{itemize}
\end{proposition}

\textit{Sketch of the proof of Proposition \ref{comparison.function}.} 
We distinguish two cases: 

\textit{Case 1:} Suppose that $z \in \big ( 0,\frac{4\pi}{3-2\varepsilon} \big ]$. We can find a smooth function $h_z: (0,\gamma(z)) \to (0,z^{\frac{3}{2}})$ such that $\lim_{s \to 0} h_z(s) = z^{\frac{3}{2}}$, $\lim_{s \to \gamma(z)} h_z(s) = 0$, and 
\begin{equation} 
\label{ode.1}
h_z'(s) = -\big ( 9\varepsilon \, z - 9\varepsilon \, h_z(s)^{\frac{2}{3}} \big )^{\frac{1}{2}} 
\end{equation}
for all $s \in (0,\gamma(z))$. Clearly, $\lim_{s \to 0} h_z'(s) = 0$. Moreover, the relation (\ref{ode.1}) implies 
\[h_z''(s) = -3\varepsilon \, h_z(s)^{-\frac{1}{3}} \leq \frac{36\pi - h_z'(s)^2}{6h_z(s)} - \frac{9}{2} \, h_z(s)^{-\frac{1}{3}}\] 
for all $s \in (0,\gamma(z))$. Therefore, the function $h$ is a solution of (\ref{ode.for.f}). Hence, if we put 
\[f_z(s) = \begin{cases} h_z(s) & \text{\rm for $s \in (0,\gamma(z))$} \\ h_z(-s) & \text{\rm for $s \in (-\gamma(z),0)$} \\ z^{\frac{3}{2}} & \text{\rm for $s = 0$,} \end{cases}\] 
then the function $f_z: (-\gamma(z),\gamma(z)) \to (0,z^{\frac{3}{2}}]$ has all the required properties.

\textit{Case 2:} Suppose that $z \in \big ( \frac{4\pi}{3-2\varepsilon},4\pi \big ]$. We can find a continuously differentiable function $h_z: (0,\gamma(z)) \to (0,z^{\frac{3}{2}})$ such that $\lim_{s \to 0} h_z(s) = z^{\frac{3}{2}}$, $\lim_{s \to \gamma(z)} h_z(s) = 0$, and 
\begin{equation} 
\label{ode.2}
h_z'(s) = \begin{cases} -\big ( 36\pi - 18(1-\varepsilon) \, y(z) \, h_z(s)^{-\frac{1}{3}} - 9 \, h_z(s)^{\frac{2}{3}} \big )^{\frac{1}{2}} & \text{\rm if $h_z(s) \geq y(z)$} \\ -\big ( 36\pi - 27(1-\varepsilon) \, y(z)^{\frac{2}{3}} - 9\varepsilon \, h_z(s)^{\frac{2}{3}} \big )^{\frac{1}{2}} & \text{\rm if $h_z(s) \leq y(z)$} \end{cases} 
\end{equation}
for all $s \in (0,\gamma(z))$. It is easy to see that $\lim_{s \to 0} h_z'(s) = 0$. For abbreviation, let $\beta(z) = \sup \{s \in (0,\gamma(z)): h_z(s) \geq y(z)\}$. Using the relation (\ref{ode.2}), we obtain 
\[h_z''(s) = \frac{36\pi - h_z'(s)^2}{6h_z(s)} - \frac{9}{2} \, h_z(s)^{-\frac{1}{3}} \leq -3\varepsilon \, h_z(s)^{-\frac{1}{3}}\] 
for all $s \in (0,\beta(z))$ and
\[h_z''(s) = -3\varepsilon \, h_z(s)^{-\frac{1}{3}} \leq \frac{36\pi - h_z'(s)^2}{6h_z(s)} - \frac{9}{2} \, h_z(s)^{-\frac{1}{3}}\] 
for all $s \in (\beta(z),\gamma(z))$. Therefore, the function $h$ is twice continuously differentiable, and satisfies the differential equation (\ref{ode.for.f}). Hence, if we define 
\[f_z(s) = \begin{cases} h_z(s) & \text{\rm for $s \in (0,\gamma(z))$} \\ h_z(-s) & \text{\rm for $s \in (-\gamma(z),0)$} \\ z^{\frac{3}{2}} & \text{\rm for $s = 0$,} \end{cases}\] 
then the function $f_z: (0,\infty) \to (0,z^{\frac{3}{2}}]$ has all the required properties. \\

\begin{lemma}
\label{asymptotics.f}
We have $\gamma(4\pi) = \pi^2$. Moreover, we have 
\[f_{4\pi}(s - \pi^2) =  6\sqrt{\pi} \, s \, \bigg [ 1 - \frac{3}{10} \, \Big ( \frac{3s}{4\pi} \Big )^{\frac{2}{3}} - \frac{3}{280} \, \Big ( \frac{3s}{4\pi} \Big )^{\frac{4}{3}} + O(s^2) \bigg ]\] 
if $s>0$ is small.
\end{lemma} 

\textit{Sketch of the proof of Lemma \ref{asymptotics.f}.} 
Using the substitution $x = (4\pi)^{\frac{3}{2}} \, \sin^3(r)$, we obtain 
\[\gamma(4\pi) = \int_0^{(4\pi)^{\frac{3}{2}}} \big ( 36\pi - 9 \, x^{\frac{2}{3}} \big )^{-\frac{1}{2}} \, dx = 4\pi \int_0^{\frac{\pi}{2}} \sin^2(r) \, dr = \pi^2.\] 
This proves the first statement. We next analyze the asymptotic behavior of the function $f_{4\pi}$ as $s \to -\pi^2$. The function $h_{4\pi}: (0,\pi^2) \to \big ( 0,(4\pi)^{\frac{3}{2}} \big )$ satisfies the differential equation 
\[h_{4\pi}'(s) = -\big ( 36\pi - 9 \, h_{4\pi}(s)^{\frac{2}{3}} \big )^{\frac{1}{2}}\] 
for all $s \in (0,\pi^2)$. This implies 
\[h_{4\pi} \big ( \pi^2 - 2\pi r + \pi \sin(2r) \big ) = (4\pi)^{\frac{3}{2}} \, \sin^3(r)\] 
for each $r \in (0,\frac{\pi}{2})$. Hence, if we put $s = 2\pi r - \pi \sin(2r)$, then we have 
\[s = \frac{4\pi}{3} \, r^3 \, \Big ( 1 - \frac{1}{5} \, r^2 + \frac{2}{105} \, r^4 + O(r^6) \Big )\] 
and 
\[h_{4\pi}(\pi^2 - s) = (4\pi)^{\frac{3}{2}} \, r^3 \, \Big ( 1 - \frac{1}{2} \, r^2 + \frac{13}{120} \, r^4 + O(r^6) \Big ).\] 
Consequently, 
\[h_{4\pi}(\pi^2 - s) =  6\sqrt{\pi} \, s \, \bigg [ 1 - \frac{3}{10} \, \Big ( \frac{3s}{4\pi} \Big )^{\frac{2}{3}} - \frac{3}{280} \, \Big ( \frac{3s}{4\pi} \Big )^{\frac{4}{3}} + O(s^2) \bigg ],\] 
as claimed. \\

\begin{lemma}
\label{bound.for.gamma}
We have $\gamma(z) < \pi^2$ for all $z \in (0,4\pi)$.
\end{lemma} 

\textit{Sketch of the proof of Lemma \ref{bound.for.gamma}.} 
The assumption (\ref{varepsilon}) implies that $\gamma(z) < \pi^2$ for all $z \in \big [ \frac{4\pi}{3-2\varepsilon},4\pi \big )$. In particular, we have $\gamma \big (  \frac{4\pi}{3-2\varepsilon} \big ) < \pi^2$. We next consider a real number $z \in (0,\frac{4\pi}{3-2\varepsilon}]$. Using the substitution $x = z^{\frac{3}{2}} \, \sin^3(r)$, we obtain 
\[\gamma(z) = \int_0^{z^{\frac{3}{2}}} \big ( 9\varepsilon \, z - 9\varepsilon \, x^{\frac{2}{3}} \big )^{-\frac{1}{2}} \, dx = \frac{1}{\sqrt{\varepsilon}} \, z \int_0^{\frac{\pi}{2}} \sin^2(r) \, dr = \frac{\pi}{4\sqrt{\varepsilon}} \, z\] 
for $z \in (0,\frac{4\pi}{3-2\varepsilon}]$. Since $\gamma \big (  \frac{4\pi}{3-2\varepsilon} \big ) < \pi^2$, we conclude that $\gamma(z) < \pi^2$ for all $z \in \big ( 0,\frac{4\pi}{3-2\varepsilon} \big ]$. This completes the proof of Lemma \ref{bound.for.gamma}. \\

We now assume that $(M,g)$ is a compact three-manifold such that $R_g \geq 6$, $\text{\rm Ric}_g \geq 2\varepsilon \, g$, and 
\[\text{\rm vol}(M,g) \geq \text{\rm vol}(S^3,\overline{g}) = 2\pi^2.\] 
Let $A: (0,2\pi^2) \to (0,\infty)$ denote the isoperimetric profile of $(M,g)$. More precisely, for each $s \in (0,2\pi^2)$, we define 
\[A(s) = \inf \big \{ \text{\rm area}(\partial \Omega,g): \text{\rm vol}(\Omega,g) = s \big \}.\] 
The following result play a key role in Bray's argument. 

\begin{proposition}[H.~Bray \cite{Bray}]
\label{ode}
Given any real number $s_0 \in (0,2\pi^2)$, there exists a real number $\delta > 0$ and a smooth function $u: (s_0-\delta,s_0+\delta) \to (0,\infty)$ with the following properties: 
\begin{itemize}
\item $u(s_0) = A(s_0)$.
\item $u(s) \geq A(s)$ for all $s \in (s_0-\delta,s_0+\delta)$.
\item $u''(s_0) \leq \min \big \{  \frac{4\pi}{u(s_0)^2} - \frac{3u'(s_0)^2}{4u(s_0)} - \frac{3}{u(s_0)}, -\frac{u'(s_0)^2}{2u(s_0)} - \frac{2\varepsilon}{u(s_0)} \big \}$.
\end{itemize}
\end{proposition}

\textit{Sketch of the proof of Proposition \ref{ode}.}  
We can find a smooth domain $\Omega \subset M$ such that $\text{\rm vol}(\Omega,g) = s_0$ and $\text{\rm area}(\partial \Omega,g) = A(s_0)$. Let $\nu$ denote the outward-pointing unit normal vector field along $\partial \Omega$. There exists a smooth one-parameter family of diffeomorphisms $\varphi_s: M \to M$, $s \in (s_0-\delta,s_0+\delta)$, with the following properties: 
\begin{itemize}
\item $\varphi_{s_0}(x) = x$ for each point $x \in M$. 
\item $\frac{\partial}{\partial s} \varphi_s(x) \big |_{s=s_0} =\frac{1}{A(s_0)} \, \nu(x)$ for each point $x \in \partial \Omega$.
\item $\text{\rm vol}(\varphi_s(\Omega),g) = s$ for each $s \in (s_0-\delta,s_0+\delta)$.
\end{itemize} 
We define a function $u: (s_0-\delta,s_0+\delta) \to (0,\infty)$ by 
\[u(s) = \text{\rm area}(\varphi_s(\partial \Omega),g).\] 
Clearly, $u(s_0) = A(s_0)$ and $u(s) \geq A(s)$ for each $s \in (s_0-\delta,s_0+\delta)$. It follows from the first variation formula that $\partial \Omega$ has constant mean curvature $H = u'(s_0)$. Using the formula for the second variation of area, we obtain 
\begin{align*} 
&\frac{d^2}{ds^2} \Big ( \text{\rm area}(\varphi_s(\partial \Omega),g) - u'(s_0) \, \text{\rm vol}(\varphi_s(\Omega),g) \Big ) \Big |_{s=s_0} \\ 
&= -\frac{1}{A(s_0)^2} \int_{\partial \Omega} (|I\!I|^2 + \text{\rm Ric}(\nu,\nu)), 
\end{align*}
where $I\!I$ denotes the second fundamental form of the boundary $\partial \Omega$. Since $\text{\rm vol}(\varphi_s(\Omega),g) = s$ for all $s \in (s_0-\delta,s_0+\delta)$, we conclude that 
\begin{equation} 
\label{2nd.variation}
u(s_0)^2 \, u''(s_0) = -\int_{\partial \Omega} (\text{\rm Ric}_g(\nu,\nu) + |I\!I|^2). 
\end{equation} 
The inequality $\text{\rm Ric}_g \geq 2\varepsilon \, g$ implies that 
\begin{align} 
\label{inequality.1}
\int_{\partial \Omega} (\text{\rm Ric}_g(\nu,\nu) + |I\!I|^2) 
&\geq \int_{\partial \Omega} \Big ( \frac{1}{2} \, H^2 + 2\varepsilon \Big ) \notag \\ 
&= \frac{1}{2} \, u(s_0) \, u'(s_0)^2 + 2\varepsilon \, u(s_0). 
\end{align} 
Since $(M,g)$ has positive Ricci curvature, the boundary $\partial \Omega$ is connected (see \cite{Bray}, p.~73). Hence, the Gauss-Bonnet theorem implies that $\int_{\partial \Omega} K \leq 4\pi$, where $K$ denotes the Gaussian curvature of $\partial \Omega$. Using the inequality $R_g \geq 6$ and the Gauss equations, we obtain 
\begin{align} 
\label{inequality.2}
\int_{\partial \Omega} (\text{\rm Ric}_g(\nu,\nu) + |I\!I|^2) \notag 
&= \frac{1}{2} \int_{\partial \Omega} (R_g - 2K + H^2 + |I\!I|^2) \notag \\ 
&\geq -4\pi + \int_{\partial \Omega} \Big ( \frac{3}{4} \, H^2 + 3 \Big ) \\ 
&= -4\pi + \frac{3}{4} \, u(s_0) \, u'(s_0)^2 + 3 \, u(s_0) \notag
\end{align} 
(cf. \cite{Bray}, p.~74). Combining (\ref{2nd.variation}), (\ref{inequality.1}), and (\ref{inequality.2}) gives 
\[u(s_0)^2 \, u''(s_0) \leq \min \Big \{ 4\pi - \frac{3}{4} \, u(s_0) \, u'(s_0)^2 - 3 \, u(s_0),-\frac{1}{2} \, u(s_0) \, u'(s_0)^2 - 2\varepsilon \, u(s_0) \Big \},\] 
as claimed. \\

Note that the function $A(s)$ may not be differentiable. For our purposes, it is sufficient that $A(s)$ is continuous.

\begin{lemma}
\label{continuity}
The function $A: (0,2\pi^2) \to (0,\infty)$ is continuous.
\end{lemma}

\textit{Sketch of the proof of Lemma \ref{continuity}.} 
Let us fix a real number $s_0 \in (0,2\pi^2)$, and let $s_k$ be a sequence of real numbers with $\lim_{k \to \infty} s_k = s_0$. It follows from Proposition \ref{ode} that $\limsup_{k \to \infty} A(s_k) \leq A(s_0)$. Hence, it suffices to show that $\liminf_{k \to \infty} A(s_k) \geq A(s_0)$. Suppose this is false. After passing to a subsequence, we may assume that the sequence $A(s_k)$ converges  to a real number $\alpha < A(s_0)$. For each $k$, we can find a smooth domain $\Omega_k \subset M$ such that $\text{\rm vol}(\Omega_k,g) = s_k$ and $\text{\rm area}(\partial \Omega_k,g) = A(s_k)$. After passing to another subsequence, we may assume that the sequence $\Omega_k$ converges to a smooth domain $\Omega_0 \subset M$ (see \cite{Ros}, Proposition 2.15). Then $\text{\rm vol}(\Omega_0,g) = s_0$ and $\text{\rm area}(\partial \Omega_0,g) = \alpha < A(s_0)$. This is a contradiction. \\

Using a comparison argument, one can relate the isoperimetric profile of $(M,g)$ to the function $f_{4\pi}$ constructed above. For abbreviation, we define a function $F: (0,2\pi^2) \to (0,\infty)$ by $F(s) = A(s)^{\frac{3}{2}}$. 

\begin{proposition}
\label{estimate.for.isoperimetric.profile}
We have $F(s) \geq f_{4\pi}(s - \pi^2)$ for all $s \in (0,2\pi^2)$.
\end{proposition}

\textit{Sketch of the proof of Proposition \ref{estimate.for.isoperimetric.profile}.} 
Suppose that $F(s_1) < f_{4\pi}(s_1 - \pi^2)$ for some real number $s_1 \in (0,2\pi^2)$. By continuity, there exists a real number $z_1 \in (0,4\pi)$ such that $s_1 \in (\pi^2 - \gamma(z_1),\pi^2 + \gamma(z_1))$ and $F(s_1) < f_{z_1}(s_1 - \pi^2)$. By Lemma \ref{bound.for.gamma}, we can find a real number $\lambda > 1$ such that $\lambda \, \gamma(z) < \pi^2$ for all $z \in (0,z_1]$ and $F(s_1) < \lambda \, f_{z_1}(\frac{s_1 - \pi^2}{\lambda})$. 

Let $I$ denote the set of all real numbers $z \in (0,z_1]$ with the property that 
\[F(s) \geq \lambda \, f_z \Big ( \frac{s - \pi^2}{\lambda} \Big )\] 
for all $s \in (\pi^2 - \lambda \, \gamma(z),\pi^2 + \lambda \, \gamma(z))$. It is easy to see that $I$ is closed and $z_1 \notin I$. Moreover, we have $(0,\alpha) \subset I$ if $\alpha > 0$ is sufficiently small. Let $z_0 \in (0,z_1)$ denote the supremum of the set $I$. Then 
\[\inf_{s \in (\pi^2 - \lambda \, \gamma(z_0),\pi^2 + \lambda \, \gamma(z_0))} \bigg [ F(s) - \lambda \, f_{z_0} \Big ( \frac{s-\pi^2}{\lambda} \Big ) \bigg ] = 0.\] 
By Lemma \ref{continuity}, the function $F: (0,2\pi^2) \to (0,\infty)$ is continuous. Hence, we can find a real number $s_0 \in  (\pi^2 - \lambda \, \gamma(z_0),\pi^2 + \lambda \, \gamma(z_0))$ such that 
\[F(s_0) = \lambda \, f_{z_0} \Big ( \frac{s_0-\pi^2}{\lambda} \Big ).\] 
By Proposition \ref{ode}, we can find a smooth function $u: (s_0-\delta,s_0+\delta) \to (0,\infty)$ such that $u(s_0) = A(s_0)$, $u(s) \geq A(s)$ for all $s \in (s_0-\delta,s_0+\delta)$, and 
\begin{equation} 
\label{ode.for.u}
u''(s_0) \leq \min \Big \{  \frac{4\pi}{u(s_0)^2} - \frac{3u'(s_0)^2}{4u(s_0)} - \frac{3}{u(s_0)}, -\frac{u'(s_0)^2}{2u(s_0)} - \frac{2\varepsilon}{u(s_0)} \Big \}. 
\end{equation}
Hence, if we define $v(s) = u(s)^{\frac{3}{2}}$, then we have 
\[v(s_0) = F(s_0) = \lambda \, f_{z_0} \Big ( \frac{s_0-\pi^2}{\lambda} \Big )\] 
and 
\[v(s) \geq F(s) \geq \lambda \, f_{z_0} \Big ( \frac{s-\pi^2}{\lambda} \Big )\] 
for all $s \in (\pi^2 - \lambda \, \gamma(z_0),\pi^2 + \lambda \, \gamma(z_0)) \cap (s_0-\delta,s_0+\delta)$. This implies 
\[v'(s_0) = f_{z_0}' \Big ( \frac{s_0-\pi^2}{\lambda} \Big )\] 
and 
\[v''(s_0) \geq \frac{1}{\lambda} \,  f_{z_0}'' \Big ( \frac{s_0-\pi^2}{\lambda} \Big ).\] 
Using the differential equation (\ref{ode.for.f}), we obtain 
\[v''(s_0) \geq \min \Big \{ \frac{36\pi - v'(s_0)^2}{6v(s_0)} - \frac{9}{2} \, \lambda^{-\frac{2}{3}} \, v(s_0)^{-\frac{1}{3}},-3\varepsilon \, \lambda^{-\frac{2}{3}} \, v(s_0)^{-\frac{1}{3}} \Big \}.\] 
On the other hand, the inequality (\ref{ode.for.u}) implies 
\[v''(s_0) \leq \min \Big \{ \frac{36\pi - v'(s_0)^2}{6v(s_0)} - \frac{9}{2} \, v(s_0)^{-\frac{1}{3}},-3\varepsilon \, v(s_0)^{-\frac{1}{3}} \Big \}.\] 
However, these inequalities are incompatible since $\lambda > 1$. The proof of Proposition \ref{estimate.for.isoperimetric.profile} is now complete. \\

Finally, we study the asymptotic behavior of the function $F(s)$ as $s \to 0$. Following an idea of Eichmair \cite{Eichmair}, we consider small geodesic balls in order to obtain an upper bound for the isoperimetric profile of $(M,g)$.

\begin{proposition}
\label{asymptotics.F}
Fix an arbitrary point $p \in M$. If $s>0$ is small, we have 
\[F(s) \leq 6\sqrt{\pi} \, s \, \bigg [ 1 - \frac{3}{2} \, c_1(p) \, \Big ( \frac{3s}{4\pi} \Big )^{\frac{2}{3}} -  \Big ( \frac{35}{24} \, c_1(p)^2 + \frac{5}{2} \, c_2(p) \Big ) \, \Big ( \frac{3s}{4\pi} \Big )^{\frac{4}{3}} + O(s^2) \bigg ].\] 
Here, the coefficients $c_1(p)$ and $c_2(p)$ are defined by $c_1(p) = \frac{1}{30} \, R_g(p)$ and $c_2(p) = \frac{1}{6300} \, (\Delta_g R_g(p) + 2 \, |\text{\rm Ric}_g(p)|^2 - 4 \, R_g(p)^2)$. 
\end{proposition}

\textit{Sketch of the proof of Proposition \ref{asymptotics.F}.} 
It follows from work of Gray and Vanhecke \cite{Gray-Vanhecke} that 
\[\text{\rm vol}(B(p,r),g) = \frac{4\pi}{3} \, r^3 \, \big ( 1 - c_1(p) \, r^2 - c_2(p) \, r^4 + O(r^6) \big )\] 
and 
\[\text{\rm area}(\partial B(p,r),g) = 4\pi \, r^2 \, \Big ( 1 - \frac{5}{3} \, c_1(p) \, r^2 - \frac{7}{3} \, c_2(p) \, r^4 + O(r^6) \Big ).\] 
Hence, if we put $s = \text{\rm vol}(B(p,r),g)$, then we obtain 
\begin{align*} 
&\text{\rm area}(\partial B(p,r),g)^{\frac{3}{2}} \\ 
&= 6\sqrt{\pi} \, s \, \bigg [ 1 - \frac{3}{2} \, c_1(p) \, \Big ( \frac{3s}{4\pi} \Big )^{\frac{2}{3}} -  \Big ( \frac{35}{24} \, c_1(p)^2 + \frac{5}{2} \, c_2(p) \Big ) \, \Big ( \frac{3s}{4\pi} \Big )^{\frac{4}{3}} + O(s^2) \bigg ]. 
\end{align*} 
Moreover, we have $F(s) \leq \text{\rm area}(\partial B(p,r),g)^{\frac{3}{2}}$ by definition of $F(s)$. Putting these facts together, the assertion follows. \\

We now complete the proof of Theorem \ref{vol.comp}. Combining Lemma \ref{asymptotics.f}, Proposition \ref{estimate.for.isoperimetric.profile}, and Proposition \ref{asymptotics.F}, we obtain 
\begin{align*} 
0 &\leq \frac{F(s) - f_{4\pi}(s - \pi^2)}{6\sqrt{\pi} \, s} \\ 
&\leq \Big ( \frac{3}{10} - \frac{3}{2} \, c_1(p) \Big ) \, \Big ( \frac{3s}{4\pi} \Big )^{\frac{2}{3}} \\ 
&+ \Big ( \frac{3}{280} - \frac{35}{24} \, c_1(p)^2 - \frac{5}{2} \, c_2(p) \Big ) \, \Big ( \frac{3s}{4\pi} \Big )^{\frac{4}{3}} + O(s^2) 
\end{align*} 
if $s>0$ is small enough. From this, we deduce that $\frac{3}{10} - \frac{3}{2} \, c_1(p) \geq 0$, hence $R_g(p) \leq 6$. On the other hand, we have $R_g(p) \geq 6$ by assumption. Since the point $p \in M$ is arbitrary, we conclude that $(M,g)$ has constant scalar curvature $6$. We next observe that $\frac{3}{280} - \frac{35}{24} \, c_1(p)^2 - \frac{5}{2} \, c_2(p) \geq 0$. This implies $c_2(p) \leq -\frac{2}{105}$, hence $|\text{\rm Ric}_g(p)|^2 \leq 12$. Since $R_g(p) = 6$, it follows that $\text{\rm Ric}_g = 2g$ for each point $p \in M$. Consequently, the manifold $(M,g)$ has constant sectional curvature $1$. Since $\text{\rm vol}(M,g) \geq \text{\rm vol}(S^3,\overline{g})$, we conclude that $(M,g)$ is isometric to $(S^3,\overline{g})$.


\begin{thebibliography}{99}
\bibitem{Andersson-Dahl}
L.~Andersson and M.~Dahl, \textit{Scalar curvature rigidity for asymptotically locally hyperbolic manifolds,} Ann. Global Anal. Geom. 16, 1--27 (1998)

\bibitem{Arnowitt-Deser-Misner}
R.~Arnowitt, S.~Deser, and C.~Misner, \textit{Coordinate invariance and energy expressions in general relativity,} Phys. Rev. 122, 997--1006 (1961)

\bibitem{Bartnik1}
R.~Bartnik, \textit{The mass of an asymptotically flat manifold,} Comm. Pure Appl. Math. 39, 661--693 (1986)

\bibitem{Bartnik2}
R.~Bartnik, \textit{Energy in general relativity,} Tsing Hua Lectures on Geometry and Analysis (Hsinchu 1990--1991), 5--27, Intl. Press, Cambridge MA (1997)

\bibitem{Besse}
A.~Besse, \textit{Einstein manifolds,} Classics in Mathematics, Springer-Verlag, Berlin (2008)

\bibitem{Boualem-Herzlich}
H.~Boualem and M.~Herzlich, \textit{Rigidity at infinity for even-dimensional asymptotically complex hyperbolic spaces,} Ann. Scuola Norm. Sup. Pisa (Ser. V), 1, 461--469 (2002)

\bibitem{Bray}
H.~Bray, \textit{The Penrose inequality in general relativity and volume comparison theorems involving scalar curvature,} PhD thesis, Stanford University (1997)

\bibitem{Bray-Brendle-Eichmair-Neves}
H.~Bray, S.~Brendle, M.~Eichmair, and A.~Neves, \textit{Area-minimizing projective planes in three-manifolds,} Comm. Pure Appl. Math. 63, 1237--1247 (2010)

\bibitem{Bray-Brendle-Neves}
H.~Bray, S.~Brendle, and A.~Neves, \textit{Rigidity of area-minimizing two-spheres in three-manifolds,} Comm. Anal. Geom. 18, 821--830 (2010)

\bibitem{Brendle}
S.~Brendle, \textit{Blow-up phenomena for the Yamabe equation,} J. Amer. Math. Soc. 21, 951--979 (2008)

\bibitem{Brendle-Marques1}
S.~Brendle and F.C.~Marques, \textit{Blow-up phenomena for the Yamabe equation II,} J. Diff. Geom. 81, 225--250 (2009)

\bibitem{Brendle-Marques2}
S.~Brendle and F.C.~Marques, \textit{Scalar curvature rigidity of geodesic balls in $S^n$,} J. Diff. Geom. 88, 379--394 (2011)

\bibitem{Brendle-Marques-Neves}
S.~Brendle, F.C.~Marques, and A.~Neves, \textit{Deformations of the hemisphere that increase scalar curvature,} Invent. Math. 185, 175--197 (2011)

\bibitem{Cai-Galloway}
M.~Cai and G.~Galloway, \textit{Rigidity of area minimizing tori in $3$-manifolds of nonnegative scalar curvature,} Comm. Anal. Geom. 8, 565--573 (2000)

\bibitem{Chrusciel-Herzlich}
P.T.~Chru\'sciel and M.~Herzlich, \textit{The mass of asymptotically hyperbolic Riemannian manifolds,} Pacific J. Math. 212, 231--264 (2003)

\bibitem{Chrusciel-Nagy}
P.T.~Chru\'sciel and G.~Nagy, \textit{The mass of spacelike hypersurfaces in asymptotically anti-de-Sitter space-times,} Adv. Theor. Math. Phys. 5, 697--754 (2001)

\bibitem{Corvino}
J.~Corvino, \textit{Scalar curvature deformation and a gluing construction for the Einstein constraint equations,} Comm. Math. Phys. 214, 137--189 (2000)

\bibitem{Eichmair}
M.~Eichmair, \textit{The size of isoperimetric surfaces in $3$-manifolds and a rigidity result for the upper hemisphere,} Proc. Amer. Math. Soc. 137, 2733--2740 (2009)

\bibitem{Fischer-Marsden}
A.E.~Fischer and J.E.~Marsden, \textit{Deformations of the scalar curvature,} Duke Math. J. 42, 519--547 (1975)

\bibitem{Fischer-Colbrie-Schoen}
D.~Fischer-Colbrie and R.~Schoen, \textit{The structure of complete stable minimal surfaces in $3$-manifolds of nonnegative scalar curvature,} Comm. Pure Appl. Math. 33, 199--211 (1980)

\bibitem{Gray-Vanhecke} 
A.~Gray and L.~Vanhecke, \textit{Riemannian geometry as determined by the volumes of small geodesic balls,} Acta Math. 142, 157--198 (1979)

\bibitem{Gromov}
M.~Gromov, \textit{Positive curvature, macroscopic dimension, spectral gaps and higher signatures,} Functional analysis on the eve of the 21st century, Vol. II (New Brunswick 1993), 1--213, Progr. Math., 132, Birkh\"auser, Boston (1996)

\bibitem{Gromov-Lawson1}
M.~Gromov and H.B.~Lawson, \textit{Spin and scalar curvature in the presence of a fundamental group,} Ann. of Math. 111, 209--230 (1980)

\bibitem{Gromov-Lawson2}
M.~Gromov and H.B.~Lawson, \textit{Positive scalar curvature and the Dirac operator on complete Riemannian manifolds,} Publ. Math. IH\'ES 58, 83--196 (1983)

\bibitem{Gursky-Viaclovsky}
M.~Gursky and J.~Viaclovsky, \textit{Volume comparison and the $\sigma_k$-Yamabe problem,} Adv. Math. 187, 447--487 (2004)

\bibitem{Hang-Wang1}
F.~Hang and X.~Wang, \textit{Rigidity and non-rigidity results on the sphere,} Comm. Anal. Geom. 14, 91--106 (2006)

\bibitem{Hang-Wang2}
F.~Hang and X.~Wang, \textit{Rigidity theorems for compact manifolds with boundary and positive Ricci curvature,} J. Geom. Anal. 19, 628--642 (2009)

\bibitem{Hersch}
J.~Hersch, \textit{Quatre propri\'et\'es isop\'erim\'etriques de membranes sph\'eriques homog\`enes,} C. R. Acad. Sci. Paris 270, 1645--1648 (1970)

\bibitem{Herzlich} 
M.~Herzlich, \textit{Scalar curvature and rigidity of odd-dimensional complex hyperbolic spaces,} Math. Ann. 312, 641--657 (1998)

\bibitem{Huang-Wu}
L.~Huang and D.~Wu, \textit{Rigidity theorems on hemispheres in non-positive space forms,} Comm. Anal. Geom. 18, 339--363 (2010)

\bibitem{Listing}
M.~Listing, \textit{Scalar curvature on compact symmetric spaces,} arxiv:1007.1832

\bibitem{Llarull}
M.~Llarull, \textit{Sharp estimates and the Dirac operator,} Math. Ann. 310, 55--71 (1998)

\bibitem{Lohkamp1}
J.~Lohkamp, \textit{Metrics of negative Ricci curvature,} Ann. of Math. 140, 655--683 (1994)

\bibitem{Lohkamp2}
J.~Lohkamp, \textit{Scalar curvature and hammocks,} Math. Ann. 313, 385--407 (1999)

\bibitem{Lohkamp3}
J.~Lohkamp, \textit{Positive scalar curvature in dim $\geq 8$,} C.R. Acad. Sci. Paris 343, 585--588 (2006)

\bibitem{Meeks-Simon-Yau}
W.~Meeks, L.~Simon, and S.T.~Yau, \textit{Embedded minimal surfaces, exotic spheres, and manifolds with positive Ricci curvature,} Ann. of Math. 116, 621--659 (1982)

\bibitem{Miao}
P.~Miao, \textit{Positive mass theorem on manifolds admitting corners along a hypersurface,} Adv. Theor. Math. Phys. 6, 1163--1182 (2002)

\bibitem{Micallef-Moraru}
M.~Micallef and V.~Moraru, \textit{Splitting of $3$-manifolds and rigidity of area-minimizing surfaces,} arxiv:1107.5346

\bibitem{Min-Oo1}
M.~Min-Oo, \textit{Scalar curvature rigidity of asymptotically hyperbolic spin manifolds,} Math. Ann. 285, 527--539 (1989)

\bibitem{Min-Oo2}
M.~Min-Oo, \textit{Scalar curvature rigidity of certain symmetric spaces,} Geometry, topology, and dynamics (Montreal, 1995), 127--137, CRM Proc. Lecture Notes vol. 15, Amer. Math. Soc., Providence RI, 1998

\bibitem{Nunes}
I.~Nunes, \textit{Rigidity of area-minimizing hyperbolic surfaces in three-manifolds,} arxiv:1103.4805

\bibitem{Parker-Taubes}
T.~Parker and C.H.~Taubes, \textit{On Witten's proof of the positive energy theorem,} Comm. Math. Phys. 84, 223--238 (1982)

\bibitem{Ros}
A.~Ros, \textit{The isoperimetric problem,} Global Theory of Minimal Surfaces (Proc. Clay Institute Summer School, 2001), 175--209, Amer. Math. Soc., Providence RI, 2005

\bibitem{Schoen-Yau1}
R.~Schoen and S.T.~Yau, \textit{On the proof of the positive mass conjecture in general relativity,} Comm. Math. Phys. 65, 45--76 (1979)

\bibitem{Schoen-Yau2}
R.~Schoen and S.T.~Yau, \textit{Existence of incompressible minimal surfaces and the topology of three dimensional manifolds of non-negative scalar curvature,} Ann. of Math. 110, 127--142 (1979)

\bibitem{Schoen-Yau3}
R.~Schoen and S.T.~Yau, \textit{On the structure of manifolds with positive scalar curvature,} Manuscripta Math. 28, 159--183 (1979)

\bibitem{Shi-Tam}
Y.~Shi and L.F.~Tam, \textit{Positive mass theorem and the boundary behaviors of compact manifolds with nonnegative scalar curvature,} J. Diff. Geom. 62, 79--125 (2002)

\bibitem{Toponogov}
V.~Toponogov, \textit{Evaluation of the length of a closed geodesic on a convex surface,} Dokl. Akad. Nauk. SSSR 124, 282--284 (1959)

\bibitem{Wang}
X.~Wang, \textit{The mass of asymptotically hyperbolic manifolds,} J. Diff. Geom. 57, 273--299 (2001)

\bibitem{Witten}
E.~Witten, \textit{A new proof of the positive energy theorem,} Comm. Math. Phys. 80, 381--402 (1981)
\end{thebibliography}
\end{document}